\documentclass[journal]{IEEEtran}
\usepackage{cite}

\ifCLASSINFOpdf
   \usepackage[pdftex]{graphicx}
\else
\fi
\usepackage{amsmath}
\usepackage{array}
\usepackage{amsfonts,setspace}
\usepackage{amsmath}
\usepackage{amssymb}
\usepackage{color}
\usepackage{float}
\usepackage{mathrsfs}
\usepackage{mathtools}
\usepackage{latexsym}
\usepackage{subfig}
\usepackage{bm}

\usepackage{enumitem}
\usepackage{booktabs}
\usepackage{algorithm}
\usepackage[noend]{algpseudocode}
\usepackage{xfrac} 
\usepackage{flexisym}
\usepackage{lineno,hyperref}
\hypersetup{
    colorlinks,
    citecolor=blue,
    filecolor=magenta,
    linkcolor=red,
}

\newcommand{\ra}{\rangle}
\newcommand{\la}{\langle}

\newcommand{\Psen}{p_{\rm p}}
\newcommand{\Psub}{p_{\rm b}}

\newcommand{\rhosen}{\rho_{\rm p}}
\newcommand{\rhosub}{\rho_{\rm b}}

\newcommand{\Csen}{c_{\rm p}}
\newcommand{\Csub}{c_{\rm b}}

\newcommand{\HH}{\mathcal{H}}
\newcommand{\FF}{\mathcal{F}}

\newtheorem{theorem}{Theorem}
\newtheorem{lemma}{Lemma}
\newtheorem{Assumption}{Assumption}

\begin{document}
%
\title{Solvability for Photoacoustic Imaging with Idealized Piezoelectric Sensors}
%
%
%

\author{Sebasti\'{a}n~Acosta~\IEEEmembership{}
\thanks{S. Acosta, Department of Pediatrics, Baylor College of Medicine and Predictive Analytics Laboratory, Texas Children's Hospital, Houston TX, USA. (sebastian.acosta@bcm.edu)}
\thanks{}}

%
%

\markboth{}%
{}

%



\maketitle

\begin{abstract}
Most reconstruction algorithms for photoacoustic imaging assume that the pressure field is measured by ultrasound sensors placed on a detection surface. However, such sensors do not measure pressure exactly due to their non-uniform directional and frequency responses, and resolution limitations. This is the case for piezoelectric sensors that are commonly employed for photoacoustic imaging. In this paper, using the method of matched asymptotic expansions and the basic constitutive relations for piezoelectricity, we propose a simple mathematical model for piezoelectric transducers. The approach simultaneously models how the pressure waves induce the piezoelectric measurements and how the presence of the sensors affects the pressure waves. Using this model, we analyze whether the data gathered by piezoelectric sensors leads to the mathematical solvability of the photoacoustic imaging problem. We conclude that this imaging problem is well-posed in certain normed spaces and under a geometric assumption. We also propose an iterative reconstruction algorithm that incorporates the model for piezoelectric measurements. A numerical implementation of the reconstruction algorithm is presented.
\end{abstract}


%
\IEEEpeerreviewmaketitle


\section{Introduction}
%
%
%
%

\IEEEPARstart{P}{hotoacoustic} tomography is a non-ionizing imaging modality designed to advantageously combine the high contrast of optical absorption with the high resolution from broadband ultrasound waves. The imaging of optical absorption reveals important functional and pathological information about biological tissues \cite{Wang2003,Beard2011,Wang-Anastasio-2011,Wang2012}. 

One of the open challenges concerning photoacoustic inversion is the incorporation of realistic models for acoustic measurements. This need for modeling the physics of ultrasound sensors has been recognized in \cite{Johansson2003,Finch2005p,Wang2007,Acosta-Montalto-2015}. It has been claimed that ultrasound measurements can be described as a linear combination of the pressure field and its normal derivative at the boundary. With that motivation, Dreier and Haltmeier \cite{Dreier2019} recently established explicit formulas for the inversion of the two-dimensional wave equation from Neumann boundary data for circular and elliptical domains. In a related effort, Zangerl, Moon and  Haltmeier \cite{Zangerl2019} derived Fourier-based reconstruction formulas for the spherical detection geometry from knowledge of Robin boundary data. 

Most other reconstruction algorithm assume that waves propagate freely across the detection boundary and that the pressure field (Dirichlet data) is measured exactly. 
These assumptions are not satisfied in practice. The pressure wave are affected by the presence of the sensors and the acoustic sensors do not measure pressure directly. They measure a surrogate for pressure that depends on the actual transducer mechanism. The most common mechanisms for ultrasound applications are based on the piezoelectric effect \cite{Wilkens2007,Wear2018b,Wear2019a}, on Fabry--Perot interferometry \cite{Beard1999b,Cox2007c,Zhang2008,Guggenheim2017} or on fiber-optic refractometry \cite{Wild2008,Wear2018a,Wissmeyer2018}. In \cite{Acosta2019a} we formulated a model specifically tailored to the Fabry--Perot sensor design. In this paper, we derive a similar model for piezoelectric sensors and analyze the well-posedness of photoacoustic imaging with such measurements. In order to attain a balance between accuracy and simplicity, the model we develop here is based on the following underlying idealizations:

\begin{enumerate}[label=(\alph*)]
\item The sensors are treated as point-like detectors. Hence, we do not account for resolution limitations due to the finite size of sensing elements. See \cite{Xu2003c,Haltmeier2010,Wang2011a,Roitner2014,Wear2018b,Wear2019a} for investigations concerning this issue. We also assume that the domain to be imaged is fully enclosed by the detection surface.
\item We assume that in each detector, the piezoelectric film is flat and its thickness is small in comparison to the wavelengths under consideration. In practice, industrial processes can manufacture piezoelectric films with thickness 30 -- 100 $\mu$m approximately \cite{Fay1994,Sirohi2000,Su2005,Wilkens2007,Oreilly2010}. 
\item Although the sensors contain elastic materials that may support shear waves, our analysis is valid for compressional waves governed by the scalar wave equation. 
\item For the piezoelectric film, the poling direction is along its thickness, and the piezoelectric properties are transversely isotropic in the plane perpendicular to the poling direction. The sensing film is mechanically isotropic.
\item Sensors may have complex structures, including a casing for structural integrity, electrodes, bonding layers and multiple paddings designed to match the mechanical impedance of the acoustic medium \cite{Brown2000,Tichy2010,Szabo2014}. However, we assume a simple design consisting of the piezoelectric film, sandwiched by electrode foils of negligible thickness, mounted on a much thicker backing layer. This follows models described in \cite{Fay1994,Sirohi2000,Szabo2014}.
\end{enumerate}

\begin{figure}[htbp]
  \centering
  \includegraphics[width=0.45\textwidth, trim = 0 -10 0 -100, clip]{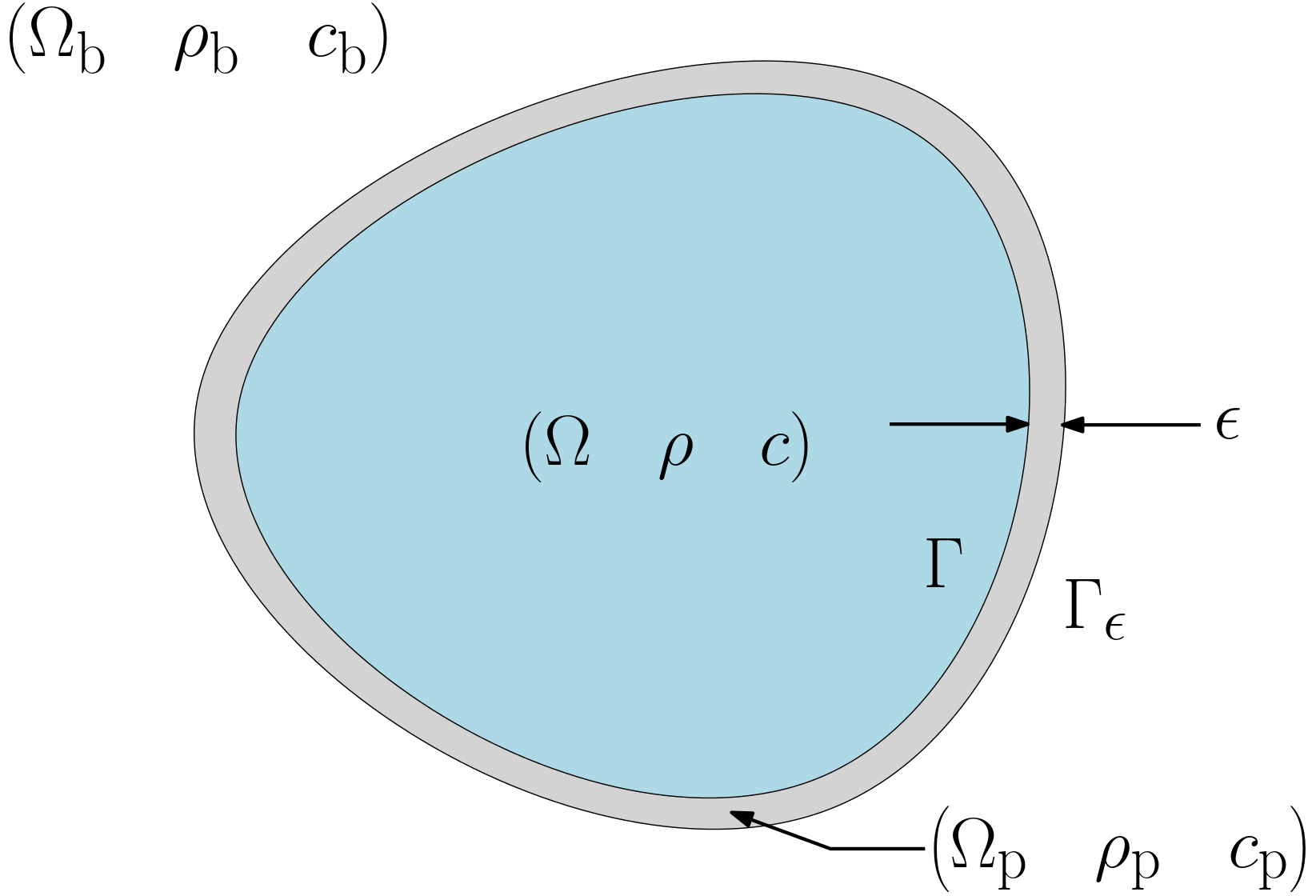}
  \caption{\textit{Acoustic domain $\Omega$ with density $\rho$ and wave speed $c$. Piezoelectric material $\Omega_{\rm p}$ of uniform thickness $\epsilon$, density $\rhosen$ and wave speed $\Csen$. The thick backing layer $\Omega_{\rm b}$ has density $\rhosub$ and wave speed $\Csub$. The interface between the acoustic domain and the piezoelectric film is denoted by $\Gamma$. The interface between the piezoelectric film and the backing layer is denoted by $\Gamma_{\epsilon}$.}}
  \label{fig:Domain}
\end{figure} 

An illustration of the idealized setup is shown in Figure \ref{fig:Domain}. The acoustic domain, denoted by $\Omega$, contains soft tissue with variable density $\rho$ and variable wave speed $c$. The piezoelectric film $\Omega_{\rm p}$ has a small uniform thickness $\epsilon > 0$, constant density $\rhosen$ and constant wave speed $\Csen$. The thick backing layer $\Omega_{\rm b}$ has constant density $\rhosub$ and constant wave speed $\Csub$. The interface between the acoustic domain and the piezoelectric film is denoted $\Gamma$. The interface between the piezoelectric film and the backing layer is denoted $\Gamma_{\epsilon}$. 

In section \ref{Section.ModelPiezo} we derive a model for the transduction from pressure to electrical voltage which is the physical quantity acquired by the piezoelectric sensors. In section \ref{Section.ModelWaves} we derive an effective boundary condition for the transmission of waves from the acoustic medium of interest into the piezoelectric sensor. This effective transmission condition accounts for the influence that the sensor exerts on the acoustic waves. Using the coupled models for piezoelectric measurements and wave propagation, in section \ref{Section.ForwardProblem} we define the forward problem  associated with photoacoustic imaging. Then in section \ref{Section.InverseProblem} we state and prove the solvability of the imaging problem with piezoelectric measurements. A reconstruction algorithm is proposed in section \ref{Section.Numerical} where some numerical simulations are presented. The conclusions follow in section \ref{Section.Conclusions}.

\section{Model for Piezoelectric Measurements} \label{Section.ModelPiezo}

We start the modeling of the piezoelectric measurements from the basic constitutive relations for both piezoelectric and mechanical variables. Since the sensing material is mechanically isotropic, the $3 \times 3$ symmetric stress tensor $\bm{\sigma}$ is related to the $3 \times 3$ symmetric strain tensor $\textbf{s}$ as follows,
\begin{equation} \label{eq:Stress_Strain}
\sigma_{ij} = \lambda \delta_{ij} \left( s_{11} + s_{22} + s_{33} \right) + 2 \mu \, s_{ij}
\end{equation}
where the direction along the thickness of the piezoelectric film is denoted as the 3-axis, and the 1-axis and 2-axis are the transverse plane. Here $\delta_{ij}$ is the Kronecker delta, and $\lambda$ and $\mu$ are the first and second Lam\'{e} coefficients. The equation of mechanical motion is
\begin{equation} \label{eq:Motion}
\rhosen \partial_{t}^2 \textbf{u} = \nabla \cdot \bm{\sigma}
\end{equation}
where $\textbf{u}$ is the $3 \times 1$ material displacement vector. For irrotational deformations, i.e. in the absence of shear stress, the above equation can be simplified in order to relate the particle displacement $\textbf{u}$ to the pressure $\Psen$ in the piezoelectric film,
\begin{equation} \label{eq:U_pressure}
\rhosen \partial_{t}^2 \textbf{u} = - \nabla \Psen
\end{equation}
where the pressure $\Psen$ is defined as 
\begin{equation} \label{eq:Pressure_Def}
 \Psen = - \left( \lambda + 2 \mu \right) \text{div} \, \textbf{u}.
\end{equation}
Combining (\ref{eq:U_pressure}) and (\ref{eq:Pressure_Def}), we find that the pressure field $\Psen$ satisfies the wave equation,
\begin{equation} \label{eq:WaveEqnPiezo}
\partial_{t}^2 \Psen = \Csen^2 \Delta \Psen 
\end{equation}
where the wave speed $\Csen$ is defined by $\Csen^2 = \left( \lambda + 2 \mu \right)/\rhosen$.

The piezoelectric transducer measures the electrical voltage $V$ across the piezoelectric film generated by the mechanical deformation due to the transmitted acoustic waves. We proceed to derive the mathematical relationship between the voltage $V$ and the pressure $\Psen$ in the piezoelectric material. Our guiding references are \cite[Ch. 5]{Tichy2010}, \cite[Ch. 5]{Szabo2014} and \cite{Sirohi2000}. 
Under small perturbations of field conditions, the linearized constitutive relation for the piezoelectric effect is the following 
\begin{equation} \label{eq:MainConstitutivePiezo}
\textbf{D} = \bm{\varepsilon} \, \textbf{E} + \textbf{d} \, \bm{\sigma}
\end{equation}
where $\textbf{D}$ is the $3 \times 1$ electric displacement vector (electric charge per area), $\textbf{E}$ is an externally applied $3 \times 1$ electric field (voltage per length) and $\bm{\varepsilon}$ is the $3 \times 3$ dielectric permittivity tensor (capacitance per length). Following \cite{Sirohi2000}, it is convenient to express the symmetric stress tensor $\bm{\sigma}$ (force per area) as a $6 \times 1$ vector and the piezoelectric tensor $\bf{d}$ (electric charge per force) as a $3 \times 6$ matrix, 
\begin{equation*} 
\textbf{d} = \left[
\begin{array}{cccccc}
0 & 0 & 0 & 0 & d_{15} & 0 \\ 
0 & 0 & 0 & d_{24} & 0 & 0 \\ 
d_{31} & d_{32} & d_{33} & 0 & 0 & 0
\end{array} \right], \quad 
\bm{\sigma} = \left[ \begin{array}{c}
\sigma_{11} \\ 
\sigma_{22} \\ 
\sigma_{33} \\ 
\sigma_{23} \\ 
\sigma_{31} \\ 
\sigma_{12}
\end{array} \right].
\end{equation*}
In the absence of an external electric field in (\ref{eq:MainConstitutivePiezo}), the normal electric displacement $D_{3}$ is given by
\begin{equation} \label{eq:SimpleConstitutivePiezo}
D_{3} = d_{31} \sigma_{11} + d_{32} \sigma_{22} + d_{33} \sigma_{33}.
\end{equation}
As assumed above, the piezoelectric tensor is transversely isotropic, which allows us to simplify the notation as follows $d_{\perp} = d_{31} = d_{32}$ and $d = d_{33}$. Combining the constitutive relations (\ref{eq:Stress_Strain}) and (\ref{eq:SimpleConstitutivePiezo}) we obtain the electric displacement in terms of the strain,
\begin{equation} \label{eq:D_S}
D_{3} = e_{\perp} s_{11} + e_{\perp} s_{22} + e s_{33}
\end{equation}
where $e_{\perp} = 2 d_{\perp} \left( \lambda + \mu \right) + d \lambda$ and $e = d \left( \lambda + 2 \mu \right) + 2 d_{\perp} \lambda$. As a consequence, using the definition of strain $\textbf{s}$ in terms of the displacement $\textbf{u}$, we obtain
\begin{equation} \label{eq:003}
D_{3} = e_{\perp} \text{div} \, \textbf{u} + \left( e - e_{\perp} \right) \partial_{n} \left( \textbf{n} \cdot \textbf{u} \right)
\end{equation}
where $\textbf{n}$ is the normal vector on $\Gamma$ and $\partial_{n}$ represents the derivative along the normal direction or 3-axis. Now we take two time-derivatives of (\ref{eq:003}) and combine with (\ref{eq:U_pressure})-(\ref{eq:Pressure_Def}) to obtain
\begin{equation} \label{eq:ElectDisp_Pressure}
\partial_{t}^2 D_{3} = - \left[ \frac{e_{\perp}}{\lambda + 2\mu} \partial_{t}^2 \Psen + \frac{\left( e - e_{\perp} \right)}{\rhosen} \partial_{n}^2 \Psen \right].
\end{equation}
Express $\Delta = \partial_{n}^2 + \Delta_{\perp}$ where $\Delta_{\perp}$ represents the surface Laplacian on the transverse plane, recall that $\Csen^{2} = (\lambda + 2\mu) / \rhosen$ and solve for $\partial_{n}^2 \Psen$ in (\ref{eq:WaveEqnPiezo}) to plug it into (\ref{eq:ElectDisp_Pressure}) to obtain
\begin{IEEEeqnarray}{ll} \label{eq:ElectDisp_Pressure2}
\partial_{t}^2 D_{3} &= - \left[ \frac{e_{\perp}}{\rhosen} \Csen^{-2} \partial_{t}^2 \Psen + \frac{\left( e - e_{\perp} \right)}{\rhosen} \left( \Csen^{-2}  \partial_{t}^2 \Psen -  \Delta_{\perp} \Psen \right) \right]  \nonumber \\
&= - \frac{e \Csen^{-2}}{\rhosen} \left[ \partial_{t}^2 \Psen - \frac{\left( e - e_{\perp} \right)}{e}  \Csen^{2}  \Delta_{\perp} \Psen   \right].
\end{IEEEeqnarray}

By definition of voltage as an electric potential, the voltage $V$ generated across the piezoelectric film by the 3-component of the electric displacement $\textbf{D}$ satisfies
\begin{equation} \label{eq:Voltage}
\partial_{n} V = \frac{D_{3}}{\varepsilon_{33}} 
\end{equation}
where $\varepsilon_{33}$ is the dielectric permittivity along the 3-axis. Integrating (\ref{eq:Voltage}) across the piezoelectric film and combining the result with
(\ref{eq:ElectDisp_Pressure2}), we obtain our model for the piezoelectric measurements
\begin{equation} \label{eq:Model_Measure}
\partial_{t}^2 V \propto \, \partial_{t}^2 \Psen - \kappa \Csen^{2} \Delta_{\perp} \Psen
\end{equation}
with vanishing initial state. Here the symbol $\propto$ denotes equality up to a multiplicative constant (which is typically estimated through experimental calibration). In (\ref{eq:Model_Measure}), it is assumed that the pressure field is constant across the piezoelectric film. This assumption is rigorously justified in the next section.

The symbol $\kappa$ appearing in the model (\ref{eq:Model_Measure}) is a unitless coefficient defined by the elastic and piezoelectric properties of the sensing film
\begin{IEEEeqnarray}{ll} \label{eq:Kappa}
\kappa = \frac{e - e_{\perp}}{e} &= \frac{2 \left( d - d_{\perp} \right) \mu}{ d \left( \lambda + 2 \mu  \right) + 2 d_{\perp} \lambda} \nonumber \\
&= \frac{\left( 1 - 2 \nu \right) \left( 1 - d_{\perp} / d \right) }{  1 - \nu \left( 1 - 2 d_{\perp}/d   \right)}
\end{IEEEeqnarray}
where we have expressed $\lambda = \rhosen \Csen^2 \nu / (1 - \nu)$ and $\mu = \rhosen \Csen^2 (1-2\nu)/(2-2\nu)$ in terms of Poisson's ratio $\nu$ to obtain the last equality. Common values for all these physical parameters are shown in Table \ref{tab:1} for polyvinylidene fluoride (PVDF) piezoelectric sensors. 

We note from (\ref{eq:Model_Measure}) that a theoretically perfect transduction from pressure to voltage would be attained if the coefficient $\kappa = 0$. However, due to the nature of the poling processes employed to manufacture these piezoelectric materials, the coefficients $d$ and $d_{\perp}$ have opposite signs and generally $|d| > |d_{\perp}|$. This implies that $1 < \left( 1 - d_{\perp}/d \right) < 2$. Hence, in order for $\kappa = 0$, the Poisson's ratio would have to be $\nu = 0.5$ which requires the piezoelectric material to be incompressible. In practice, Poisson's ratio for PVDF films ranges from 0.2 to 0.4 approximately. We note that $\kappa$ ranges from 0.3 to 1.5, for the realistic range of values for the Poisson's ratio $\nu$ and the piezoelectric ratio $d_{\perp}/d$ displayed in Table \ref{tab:1}.

\begin{table}[htbp]
\caption{\textit{Estimates for the physical parameters of PVDF piezoelectric sensors \cite{Fay1994,Sirohi2000,Brown2000,Zeqiri2007,Wilkens2007,Bedard2008,Szabo2014,Cao2018}.}}
\label{tab:1}       
\begin{tabular}{lrl}
\hline\noalign{\smallskip}
\textbf{Parameter} & \textbf{Value} & \textbf{Units}  \\
\noalign{\smallskip}\hline\noalign{\smallskip}
PVDF thickness $\epsilon$ 			& 10 -- 60 		& $\mu$m 		 \\
PVDF density $\rhosen$ 				& 1780 -- 1950	& kg m$^{-3}$ 	 \\
PVDF wave speed $\Csen$ 			& 1300 -- 2300	& m s$^{-1}$ 	 \\
PVDF Poisson's ratio $\nu$ 			& 0.2 -- 0.4	&  	 \\
Piezoelectric coeff. $d$ 			& -(30 -- 35)	& pC/N 	 \\
Piezoelectric coeff. $d_{\perp}$ 	& 3 -- 15		& pC/N 	 \\
Coefficient $\kappa$				& 0.3 -- 1.5	& \\ 
Backing density $\rhosub$ 			& 1900 -- 2500	& kg m$^{-3}$ 	 \\
Backing wave speed $\Csub$ 			& 1000 -- 4000	& m s$^{-1}$ 	 \\
\noalign{\smallskip}\hline
\end{tabular} 
\end{table}

\section{Effective Model for Wave Propagation} \label{Section.ModelWaves}

Typically, the piezoelectric film and the backing layer are acoustically more rigid and heavier than the biological medium of interest. Therefore, the presence of the sensors induces partial reflections of the waves. Here we seek to model how the sensors exert influence on the pressure waves. This model takes the form of an effective impedance boundary condition that replaces the more involved transmission process for waves traveling from the acoustic domain $\Omega$, through the piezoelectric film $\Omega_{\rm p}$ and into the backing layer $\Omega_{\rm b}$. We assume that the pressure field $\Psub$ in the backing layer is outgoing which translates into satisfying a radiation condition of the form,
\begin{equation} \label{eq:RadCond}
\partial_{n} \Psub + \Csub^{-1} \partial_{t} \Psub + \HH \Psub = 0 \qquad \text{on $\Gamma_{\epsilon}$}
\end{equation}
where $\HH$ is the mean curvature of the surface $\Gamma$. See \cite{Antoine1999,Acosta2017c} for a derivation.

As in \cite{Acosta2019a}, we make some geometric assumptions about the domain $\Omega_{\rm p}$ occupied by the  piezoelectric film. We let $\Omega_{\rm p} = \{ \textbf{y} \in \Omega^{\rm c} ~:~ 0 < \text{dist}(\textbf{y},\Gamma) < \epsilon \}$. For sufficiently small $\epsilon$, the domain $\Omega_{\rm p}$ can be expressed as a family of parallel surfaces parametrized by $0 < z < \epsilon$ and defined by $\Gamma_{z} = \{ \textbf{y} = \textbf{x} + z \textbf{n}(\textbf{x}) ~:~ \textbf{x} \in \Gamma \}$ where $\textbf{n}(\textbf{x})$ is the normal vector at $\textbf{x} \in \Gamma$. For smooth $\Gamma$ and sufficiently small $\epsilon$, the surfaces $\Gamma_{z}$ are well-defined, smooth and mutually disjoint. Each point $\textbf{y} \in \Omega_{\rm p}$ can be uniquely represented in the form $\textbf{y} = \textbf{x} + z \textbf{n}(\textbf{x})$ for $\textbf{x} \in \Gamma$ and $0<z<\epsilon$. In addition, the normal vector at $\textbf{y} \in \Gamma_{z}$ coincides with the normal vector at $\textbf{x} \in \Gamma$. See details concerning parallel surfaces in \cite[Sect. 6.2]{Kress-Book-1999}.

The transmission of the pressure field from the acoustic domain into the piezoelectric film is governed by the following transmission conditions at the interface $\Gamma$,
\begin{equation} \label{eq:Transmission1}
p = \Psen \quad \text{and} \quad \rho^{-1} \partial_{n} p  = \rhosen^{-1} \partial_{n} \Psen  \qquad \text{on $\Gamma$,}
\end{equation}
where $p$ and $\Psen$ are the pressure in the acoustic medium and piezoelectric film, respectively. The first condition in (\ref{eq:Transmission1}) ensures the continuity of the pressure field. The second condition in (\ref{eq:Transmission1}) ensures the continuity of particle motion in the normal direction. Similar transmission conditions hold at the interface $\Gamma_{\epsilon}$,
\begin{equation} \label{eq:Transmission2}
\Psen = \Psub \quad \text{and} \quad \rhosen^{-1} \partial_{n} \Psen = \rhosub^{-1} \partial_{n} \Psub   \qquad \text{on $\Gamma_{\epsilon}$,}
\end{equation}
where $\Psub$ is the pressure in the backing layer. The pressures $p$, $\Psen$ and $\Psub$ satisfy the wave equation with respective wave speeds $c$, $\Csen$ and $\Csub$.

Now we proceed to use the method of matched asymptotic expansions to derive an effective model for the interplay between the pressure fields and the piezoelectric sensor. For an introduction to this method, see \cite{Hinch2012}. First we consider the formal asymptotic expansions for the pressure fields,
\begin{subequations} \label{eq:Expansion}
\begin{IEEEeqnarray}{l} 
p(t,\textbf{x}) = p^{0}(t,\textbf{x}) + \epsilon \, p^{1}(t,\textbf{x}) + \mathcal{O}(\epsilon^2) \\
\Psen(t,\textbf{x},z) = \Psen^{0}(t,\textbf{x},z) + \epsilon \, \Psen^{1}(t,\textbf{x},z) + \mathcal{O}(\epsilon^2) \\
\Psub(t,\textbf{x}) = \Psub^{0}(t,\textbf{x}) + \epsilon \, \Psub^{1}(t,\textbf{x}) + \mathcal{O}(\epsilon^2)
\end{IEEEeqnarray}
\end{subequations}
and introduce a change of variable in order to extract the effect of the piezoelectric film thickness $\epsilon$,
\begin{equation} \label{eq:ChangeVariable}
z = \epsilon \, \zeta \qquad \text{for $\zeta \in [0,1]$}.
\end{equation}

The boundary value problem for the pressure field $\Psen$ in the piezoelectric film, governed by the wave equation (\ref{eq:WaveEqnPiezo}) and the transmission conditions (\ref{eq:Transmission1})-(\ref{eq:Transmission2}), can be recast in terms of $\zeta$ and terms with same powers of $\epsilon$ are gathered to obtain the following cases.
\paragraph{$\mathcal{O}(\epsilon^{0})$-terms}
\begin{IEEEeqnarray*}{l} 
\partial_{\zeta}^2 \Psen^{0} = 0, \qquad \text{for $\zeta \in (0,1)$,} \\
p^{0} = \Psen^{0} \quad \text{and} \quad \partial_{\zeta} \Psen^{0} = 0, \quad \text{at $\zeta = 0$,} \\
\Psen^{0} = \Psub^{0} \quad \text{and} \quad \partial_{\zeta} \Psen^{0} = 0, \quad \text{at $\zeta = 1$,}
\end{IEEEeqnarray*}
which imply that $\Psen^{0}$ is constant as a function of $\zeta$ and that the first effective transmission condition is that  
\begin{eqnarray} 
p^{0} = \Psen^{0} = \Psub^{0}, \qquad  \text{on $\Gamma$}. \label{eq:O-0-trans1}
\end{eqnarray}

\paragraph{$\mathcal{O}(\epsilon^{1})$-terms}
\begin{IEEEeqnarray*}{l}
\partial_{\zeta}^2 \Psen^{1} = 0, \qquad \text{for $\zeta \in (0,1)$,} \\
p^{1} = \Psen^{1} \quad \text{and} \quad \rho^{-1} \partial_{n} p^{0} = \rhosen^{-1} \partial_{\zeta} \Psen^{1} , \quad \text{at $\zeta = 0$,} \\
\Psen^{1} = \Psub^{1} \quad \text{and} \quad \rhosen^{-1} \partial_{\zeta} \Psen^{1} = \rhosub^{-1} \partial_{n} \Psub^{0}, \quad \text{at $\zeta = 1$,}
\end{IEEEeqnarray*}
which imply that $\partial_{\zeta} \Psen^{1}$ is constant as a function of $\zeta$ and that the second effective transmission condition is 
\begin{eqnarray} \label{eq:O-0-trans2}
\rho^{-1} \partial_{n} p^{0} = \rhosen^{-1} \partial_{\zeta} \Psen^{1} = \rhosub^{-1} \partial_{n} \Psub^{0} \qquad \text{on $\Gamma$}.
\end{eqnarray}

Combining (\ref{eq:RadCond}) and (\ref{eq:O-0-trans1})-(\ref{eq:O-0-trans2}), we obtain closed-form effective governing equations for the leading order term $p^{0}$ of the acoustic field in the domain $\Omega$,
\begin{IEEEeqnarray}{ll}
\partial_{t}^2 p^{0} - c^2 \Delta p^{0} = 0 &\quad \text{in $\{ t > 0 \} \times \Omega$}, \label{eq:EffBVP-01} \\
\rhosub \partial_{n} p^{0} + \rho \Csub^{-1} \partial_{t} p^{0} + \rho \HH p^{0} = 0 &\quad \text{on $\{ t > 0 \} \times \Gamma$}. \label{eq:EffBVP-02}
\end{IEEEeqnarray}
Similar models for photoacoustics are studied in \cite{Acosta-Montalto-2015,Acosta-Montalto-2016,Acosta2019a,Acosta2019b}.

\begin{figure}[h]
  \centering
  \includegraphics[width=0.45\textwidth, trim=0 15 0 30, clip]{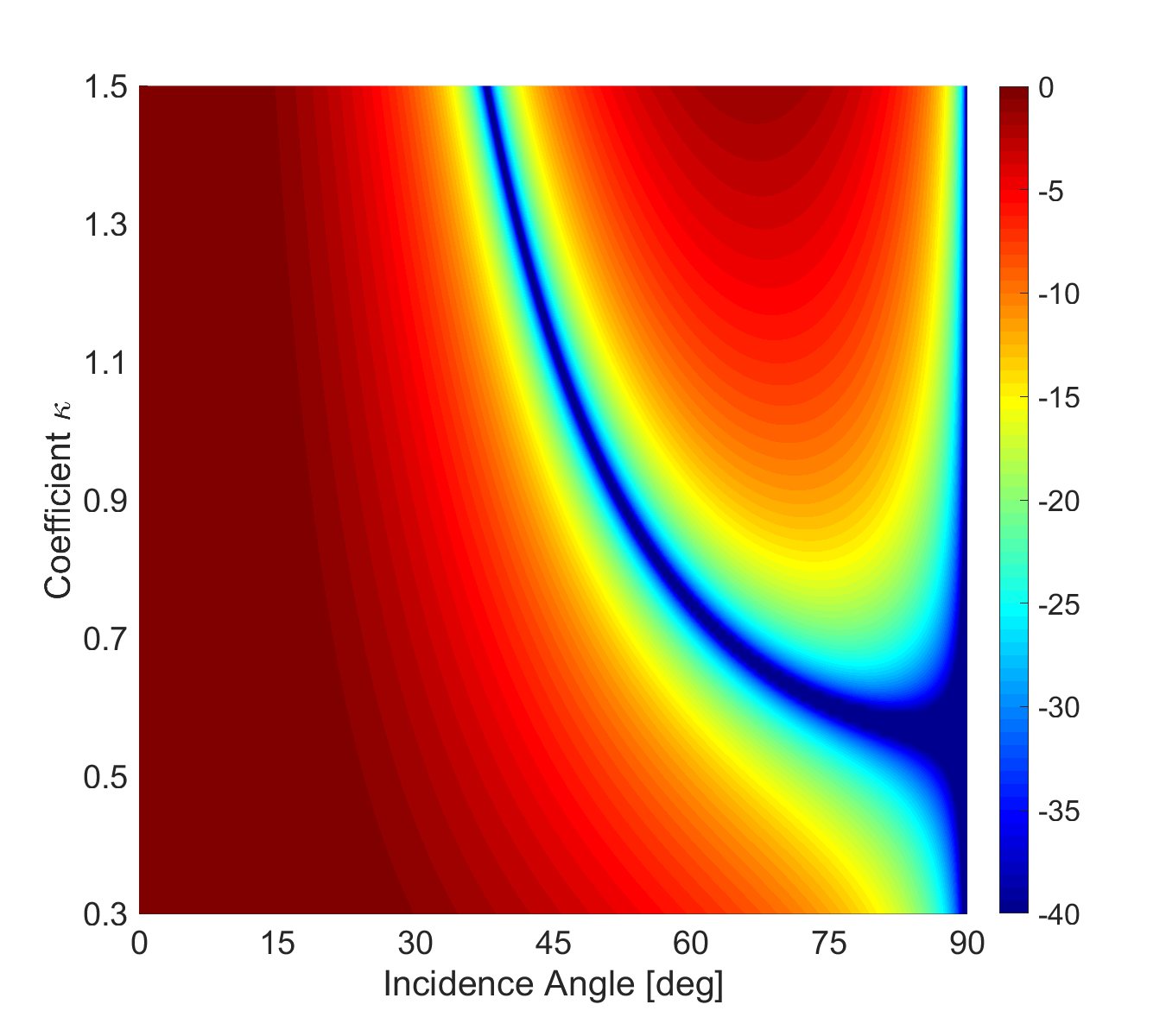}
  \caption{\textit{Directivity (\ref{eq:response}) in decibels for piezoelectric sensor as a function of the incidence angle $\theta$ and coefficient $\kappa$. The parameters correspond to a PVDF film with (compressional) wave speed $\Csen = 2000$ m/s and density $\rhosen = 1800$ kg/m$^3$ and a backing layer with
(compressional) wave speed $\Csub = 1000$ m/s and density $\rhosub = 2000$ kg/m$^3$. The acoustic medium corresponds to water with wave speed $c = 1500$ m/s and density $\rho = 1000$ kg/m$^3$.}}
  \label{fig:directivity}
\end{figure}

As an example for the response of the piezoelectric sensor design, we can analyze its behavior for plane waves and for a flat boundary $\Gamma$. Both the boundary value problem (\ref{eq:EffBVP-01})-(\ref{eq:EffBVP-02}) and the model for the measurements (\ref{eq:Model_Measure}) play an important role in this analysis. A plane wave of the form $p_{\rm inc} = e^{i(\textbf{x} \cdot \textbf{k} - \omega t)}$ propagating in the direction of $\textbf{k}$, induces a reflection governed by (\ref{eq:EffBVP-02}). The total pressure field $p$ is the superposition of the incident and reflected wave,
\begin{equation} \label{eq:plane_wave_01}
p(\textbf{x},t) = e^{i(\textbf{x} \cdot \textbf{k} - \omega t)} + R e^{i(\textbf{x} \cdot \textbf{k}_{\rm r} - \omega t)} + \mathcal{O}(\epsilon)
\end{equation}
where $R$ is the reflection coefficient, $\textbf{k}_{\rm r}$ is the reflection wavenumber satisfying $|\textbf{k}| = |\textbf{k}_{\rm r}| = \omega/c$ and $\textbf{n} \cdot \textbf{k}_{\rm r} = - \textbf{n} \cdot \textbf{k}$, where $\textbf{n}$ is the outward normal on $\Gamma$. We can write $\textbf{n} \cdot \textbf{k}= |\textbf{k}| \cos \theta$ where $\theta$ is the angle of incidence. Plugging (\ref{eq:plane_wave_01}) into (\ref{eq:EffBVP-02}) and neglecting the $\mathcal{O}(\epsilon)$ terms, we find that the reflection coefficient satisfies 
\begin{equation} \label{eq:reflection_coeff}
R = \frac{\cos \theta - \alpha}{\cos \theta + \alpha}, \qquad \text{where} \quad \alpha = \frac{\rho c}{ \rhosub \Csub}.
\end{equation}
After plugging (\ref{eq:plane_wave_01})-(\ref{eq:reflection_coeff}) into the model (\ref{eq:Model_Measure}) and evaluating at the origin $\textbf{x} = \textbf{0}$, we find that the piezoelectric measurements satisfy the following directivity pattern
\begin{equation} \label{eq:response}
\frac{V}{p_{\rm inc}} =  \left(1 + \frac{\cos \theta - \alpha}{\cos \theta + \alpha} \right) \left( 1 - \kappa \frac{\Csen^2}{c^2} \sin^{2} \theta \right)
\end{equation}
where $\kappa$ is given by (\ref{eq:Kappa}). Figure \ref{fig:directivity} displays the directional response (\ref{eq:response}) in decibels as a function of the incidence angle $\theta$ and piezoelectric coefficient $\kappa$ over a realistic range of values shown in Table \ref{tab:1}. We observe that for values of $\kappa > c^2 / \Csen^2$, a critical angle appears. For incidence at this critical angle, vanishing measurements are obtained by the piezoelectric sensor design. This critical angle is given by $\theta_{\rm cr} = \arcsin \left( c \, \Csen^{-1} \kappa^{-1/2} \right)$.

\section{The Forward Problem} \label{Section.ForwardProblem}

Now we proceed to define the forward problem of photoacoustic imaging in terms of the wave propagation model (\ref{eq:EffBVP-01})-(\ref{eq:EffBVP-02}) and the model for piezoelectric measurements (\ref{eq:Model_Measure}). We neglect higher order terms $\mathcal{O}(\epsilon)$ studied in the previous section, so that the pressure field $p$ is assumed to satisfy the following initial value problem,
\begin{subequations} \label{eq:IVBP}
\begin{IEEEeqnarray}{ll}
 \partial_{t}^2 p - c^2 \Delta p = 0 &\quad \text{in $(0,T) \times \Omega$} \label{eq:IVBP-01} \\
 \rhosub \partial_{n} p + \rho \Csub^{-1} \partial_{t} p + \rho \HH p = 0 &\quad \text{on $(0,T) \times \Gamma$} \label{eq:IVBP-02} \\
 p = f \quad \text{and} \quad \partial_{t} p = 0 &\quad \text{on $\{ t=0 \}  \times \Omega$} \label{eq:IVBP-03}
\end{IEEEeqnarray}
\end{subequations}
where $0 < T < \infty$ is the measurement time to be determined later. Recall that the underlying assumption concerning media properties are that $c$ is bounded from below and above, and is smooth in $\Omega$, and that $\Csen$, $\Csub$, $\rhosen$, $\rhosub$, $\rho$ and $\kappa$ are constants. The forward mapping is given by
\begin{eqnarray} \label{eq:forward}
\FF : f \mapsto V
\end{eqnarray}
where, according to the piezoelectric model (\ref{eq:Model_Measure}), the measured electric voltage $V$ satisfies
\begin{subequations} \label{eq:IVBP-V}
\begin{IEEEeqnarray}{ll}
\partial_{t}^2 V = \partial_{t}^2 p - \kappa \Csen^2 \Delta_{\perp} p &\qquad \text{on $(0,T) \times \Gamma$} \label{eq:IVBP-05} \\
V = \partial_{t} V = 0 &\qquad \text{on $\{ t=0 \} \times \Gamma$} \label{eq:IVBP-06}
\end{IEEEeqnarray}
\end{subequations}
for the pressure field $p$ evolving according to (\ref{eq:IVBP}) from the initial condition $f$. The mapping (\ref{eq:forward}), which we seek to invert for photoacoustic imaging, encodes the physical principles for acoustic wave propagation from the unknown pressure profile $f$ to the measured electrical voltage $V$ generated by the piezoelectric sensors.

We work with the standard Sobolev spaces based on square-integrable functions defined on the domain $\Omega$ or the boundary $(0,T) \times \Gamma$. The associated inner product extends as the duality pairing between functionals and functions. For the Sobolev space $H^{0}(\Omega)$, the inner product is weighted by $c^{-2}$ so that the differential operator $c^2 \Delta$ is formally self-adjoint. The well-posedness in Sobolev spaces of the initial value problem (\ref{eq:IVBP}) is a well-established result \cite{EvansPDE,Lio-Mag-Book-1972}.

\section{The Inverse Problem} \label{Section.InverseProblem}

The inverse problem associated with photoacoustic imaging is the following: Given the voltage measurements $V$ modeled by (\ref{eq:IVBP-V}) on $\Gamma \times (0,T)$, induced by the pressure field $p$ satisfying (\ref{eq:IVBP}), find the unknown initial condition $f$. The solvability of this inverse problem depends on the geometry of the domain $\Omega$, the profile of the wave speed $c$ and the time $T < \infty$. These conditions are made precise in the following assumption, known as the geometric control condition or a nontrapping condition for the geodesic flow. We work with the manifold $\Omega$ endowed with the Riemannian metric $c^{-2} dx^2$. See \cite{Bardos1992,GlowinskiLionsHe2008} for details.

\begin{Assumption}[Nontrapping condition] \label{Ass:Nontrapping}
Let $\Omega$ be a simply connected bounded domain with smooth boundary $\Gamma$. Assume there exists $T_{\rm o} < \infty$ such that any (unit speed) geodesic ray of the manifold $(\Omega,c^{-2} dx^2)$, originating from any point in $\Omega$ at time $t=0$, reaches the boundary $\Gamma$ at a nondiffractive point before $t=T_{\rm o}$.
\end{Assumption}

With this assumption in place, we can state the main result of the paper in the form of a theorem.
\begin{theorem} \label{Thm:Main}
Under the Assumption \ref{Ass:Nontrapping} for the manifold $(\Omega, c^{-2} dx^2)$ and time $T>T_{\rm o}$, the forward mapping $\FF : H^{1}_{0}(\Omega) \to H^{1}([0,T];H^{0}(\Gamma))$ is injective, that is, the photoacoustic imaging problem is uniquely solvable. Moreover, the following stability estimate,
\begin{equation} \label{eq:stability_PAT}
\| f \|_{H^{0}(\Omega)} \leq C \| V \|_{H^{1}([0,T];H^{0}(\Gamma))}
\end{equation}
holds for some constant $C>0$.
\end{theorem}

We wish to make some comments before we proceed with the proof. Notice in (\ref{eq:stability_PAT}) that we are only able to dominate $f$ in the norm of $H^{0}(\Omega)$ (rather than in the norm of its stated space $H^{1}_{0}(\Omega)$) with the measured data $V$ in the norm of $H^{1}([0,T];H^{0}(\Gamma))$. By contrast, when the Dirichlet data is measured on $[0,T] \times \Gamma$, then the imaging operator (left inverse of $\FF$) enjoys stability estimates as a mapping from $H^{0}([0,T] \times \Gamma)$ to $H^{0}(\Omega)$, or from $H^{1}([0,T] \times \Gamma)$ to $H^{1}_{0}(\Omega)$.
See \cite{Bardos1992,Stefanov2009,Acosta-Montalto-2015} for details. Hence, there is an apparent loss of stability due to the nature of the piezoelectric measurement model (\ref{eq:IVBP-V}). The double time-integration needed to invert the left-hand side of (\ref{eq:IVBP-05}) does not fully restore the regularity lost by the application of the hyperbolic differential operator on the right-hand side of (\ref{eq:IVBP-05}). This is a well-known property concerning regularity of hyperbolic equations \cite{EvansPDE}. We also note that Theorem \ref{Thm:Main} is slightly different from what is presented in \cite{Acosta2019a} where the imaging operator was shown to satisfy a stability estimate as a mapping from $H^{0}([0,T];H^{1}(\Gamma))$ to $H^{0}(\Omega)$. Hence, formally, there is a mild loss of stability of one degree either in space or in time, but not both.

Now, it is convenient to define the following operation
\begin{equation}
(\partial_{t}^{-1} v)(t) = \int_{0}^{t} v(s) \, ds  \label{eq:inv_partial_t}
\end{equation}
so that $\partial_{t}^{-1} \partial_{t} v = \partial_{t} \partial_{t}^{-1} v = v$ for any sufficiently smooth $v$ such that $v=0$ at $t=0$. Now let
\begin{equation}
u = \partial_{t}^{-1} p  \label{eq:u_int_p}
\end{equation}
where $p$ and $V$ satisfy (\ref{eq:IVBP-V}) and $p(0) = f \in H^{1}_{0}(\Omega)$. Then it follows that $u$ solves the following initial value problem,
\begin{subequations} \label{eq:IVBP_u}
\begin{IEEEeqnarray}{ll}
\partial_{t}^2 u - \kappa \Csen^2 \Delta_{\perp} u = \partial_{t} V &\quad \text{on $(0,T) \times \Gamma$} \\
u = \partial_{t} u = 0 &\quad \text{on $\{ t=0 \} \times \Gamma$}.
\end{IEEEeqnarray}
\end{subequations}
The following lemma is a well-established result. See \cite[\S 7.2, Thms. 3-5]{EvansPDE} or \cite[Ch. 3, \S 8, Thm. 8.1]{Lio-Mag-Book-1972} for details.
\begin{lemma} \label{Lemma:stability_u}
Let $u$ solve (\ref{eq:IVBP_u}). If $V \in H^{1}([0,T]; H^{0}(\Gamma))$, then the field $u \in C^{k}([0,T] ; H^{1-k}(\Gamma))$ for $k=0,1$. Moreover, the following stability estimate 
\begin{equation} \label{eq:stability_u}
\| u \|_{C^{k}([0,T]; H^{1-k}(\Gamma))} \leq C \| V \|_{H^{1}([0,T]; H^{0}(\Gamma))}
\end{equation}
holds for some constant $C>0$.
\end{lemma}

The definition of the Bochner spaces $H^{1}([0,T]; H^{0}(\Gamma))$ and $C^{k}([0,T]; H^{1-k}(\Gamma))$ can be found in \cite{EvansPDE,Lio-Mag-Book-1972}. Using Lemma \ref{Lemma:stability_u}, we proceed to prove the main theoretical result of the paper. In what follows, the generic constant $C>0$ changes from inequality to inequality, but it does not depend on $f$, $p$ or $V$.

\begin{IEEEproof}[Proof of Theorem \ref{Thm:Main}]
Under Assumption \ref{Ass:Nontrapping} for the manifold $(\Omega, c^{-2} dx^2)$ and time $T>T_{\rm o}$, observability of waves from the boundary \cite{Bardos1992,GlowinskiLionsHe2008} yields that
\begin{equation} \label{eq:observability}
\| f \|_{H^{0}(\Omega)} \leq C \| p \|_{H^{0}([0,T] \times \Gamma)}
\end{equation}
for some constant $C=C(\Omega,c,T)$. Now, from the definition (\ref{eq:u_int_p}) of $u$ and the stability estimate in Lemma \ref{Lemma:stability_u} for $k=1$, we obtain that
\begin{equation} \label{eq:stability_p_bdy}
\| p \|_{C^{1}([0,T]; H^{0}(\Gamma))} \leq C \| V \|_{H^{1}([0,T]; H^{0}(\Gamma))}.
\end{equation}
Since the norm of $C^{1}([0,T]; H^{0}(\Gamma))$ dominates the norm of $H^{0}([0,T] \times \Gamma)$, combining (\ref{eq:observability}) and (\ref{eq:stability_p_bdy}) we obtain the desired result (\ref{eq:stability_PAT}).
\end{IEEEproof}

\section{Numerical Simulations} \label{Section.Numerical}

Now we propose and numerically implement a reconstruction algorithm to solve
the PAT problem at the discrete level. The reconstructions presented here are based on the Landweber iterative method \cite[Ch. 6]{EnglBook2000} with Nesterov's acceleration. See Stefanov and Yang \cite{Stefanov2017a} for an excellent analysis of the Landweber method for a similar problem. The iterative process is defined in Algorithm \ref{alg:landweber}. In brief, the Landweber method is based on inverting (\ref{eq:forward}) by solving 
\begin{equation*}
\left( I - \left( I - \gamma \FF^{*} \FF \right) \right) f = \gamma \FF^{*} V
\end{equation*}
where the normal operator $\left( \FF^{*} \FF \right)$ is positive definite and $\gamma>0$ is chosen small enough so that the spectrum of $\left( I - \gamma \FF^{*} \FF \right)$ is contained in $(-1,1)$ making it a contraction. The parameter $\gamma$ is known as the relaxation factor. The parameter $\mu$, known as the momentum factor, allows for the acceleration of the convergence. Unfortunately, it is hard to choose $\gamma$ and $\mu$ optimally. We resort to trial and error to set them satisfactorily. In the presence of noise, the number $K$ of iterations is chosen according to a regularization rule.

For the Landweber method, it is necessary to evaluate the adjoint $\FF^{*}$ of the forward operator $\FF$. This evaluation amounts to solve the following final boundary value problem,
\begin{subequations} \label{eq:IVBP_adj}
\begin{IEEEeqnarray}{ll}
\partial_{t}^2 \varphi - c^2 \Delta \varphi = 0 &\quad \text{$(0,T) \times \Omega$} \quad \label{eq:IVBP-01_adj} \\
\rhosub \partial_{n} \varphi - \rho \Csub^{-1} \partial_{t} \varphi + \rho \HH \varphi =  - \rhosub \eta &\quad \text{$(0,T) \times \Gamma$} \quad \label{eq:IVBP-02_adj} \\
\varphi = 0 \quad \text{and} \quad \partial_{t} \varphi = 0 &\quad \text{$\{ t=T \}  \times \Omega$} \quad \label{eq:IVBP-03_adj}
\end{IEEEeqnarray}
\end{subequations}
where $\eta$ solves
\begin{subequations} \label{eq:IVBP_adj_eta}
\begin{IEEEeqnarray}{ll}
\partial_{t}^2 \eta = \partial_{t}^2 \psi - \kappa \Csen^2 \Delta_{\perp} \psi &\quad \text{in $(0,T) \times \Gamma$} \label{eq:IVBP-01_adj_eta} \\
\eta = 0 \quad \text{and} \quad \partial_{t} \eta = 0 &\quad \text{on $\{ t=T \}  \times \Gamma$} \label{eq:IVBP-02_adj_eta}
\end{IEEEeqnarray}
\end{subequations}
in order to define the adjoint mapping as
\begin{equation} \label{eq:adjoint}
\FF^{*} : \psi \mapsto  \partial_{t} \varphi |_{t=0}. 
\end{equation}
It can be shown, using straight-forward integration by parts, that indeed $\FF^{*}$ is the adjoint of $\FF$ or equivalently that
\begin{IEEEeqnarray*}{ll}
\la f , \partial_{t} \varphi|_{t=0} \ra_{\Omega} &= \la p , \eta \ra_{[0,T] \times \Gamma} \\
&= \la p , \left( I - \kappa \Csen^2 \Delta_{\perp} \partial_{t}^{-2} \right) \psi \ra_{[0,T] \times \Gamma} \\
&= \la V , \psi \ra_{[0,T] \times \Gamma}
\end{IEEEeqnarray*}
for all sufficiently smooth $f$ and $\psi$, where $V = \FF(f)$ according to (\ref{eq:forward}) and $\partial_{t} \varphi|_{t=0} = \FF^{*}(\psi)$ according to (\ref{eq:adjoint}). The brackets $\la \cdot , \cdot \ra_{X}$ denote the inner product of the Hilbert space $H^{0}(X)$ which extends as the duality pairing between the Sobolev functionals and functions.

\begin{algorithm}[h]
\caption{Accelerated Landweber iteration}
\label{alg:landweber}
\begin{algorithmic}
\State{Set $0 < K$, $0 < \gamma < 2 \| \FF \|^{-2}$ and $0 \leq \mu < 1$.}
\State{Initial guesses $u_{0} = \FF^{*}V$ and $v_{0} = u_{0}$}
\For{$k=1,2,..., K$}
\State{$v_{k} = u_{k-1} - \gamma \left( \FF^{*} \FF u_{k-1} - u_{0} \right) / \| u_{0} \| $}
\State{$u_{k} = v_{k} + \mu \left( v_{k} - v_{k-1}  \right)$}
\EndFor
\State{return $u_{K}$}
\end{algorithmic}
\end{algorithm}

\begin{figure}[htbp]
  \centering
  \includegraphics[width=0.45\textwidth, trim=10 30 10 45, clip]{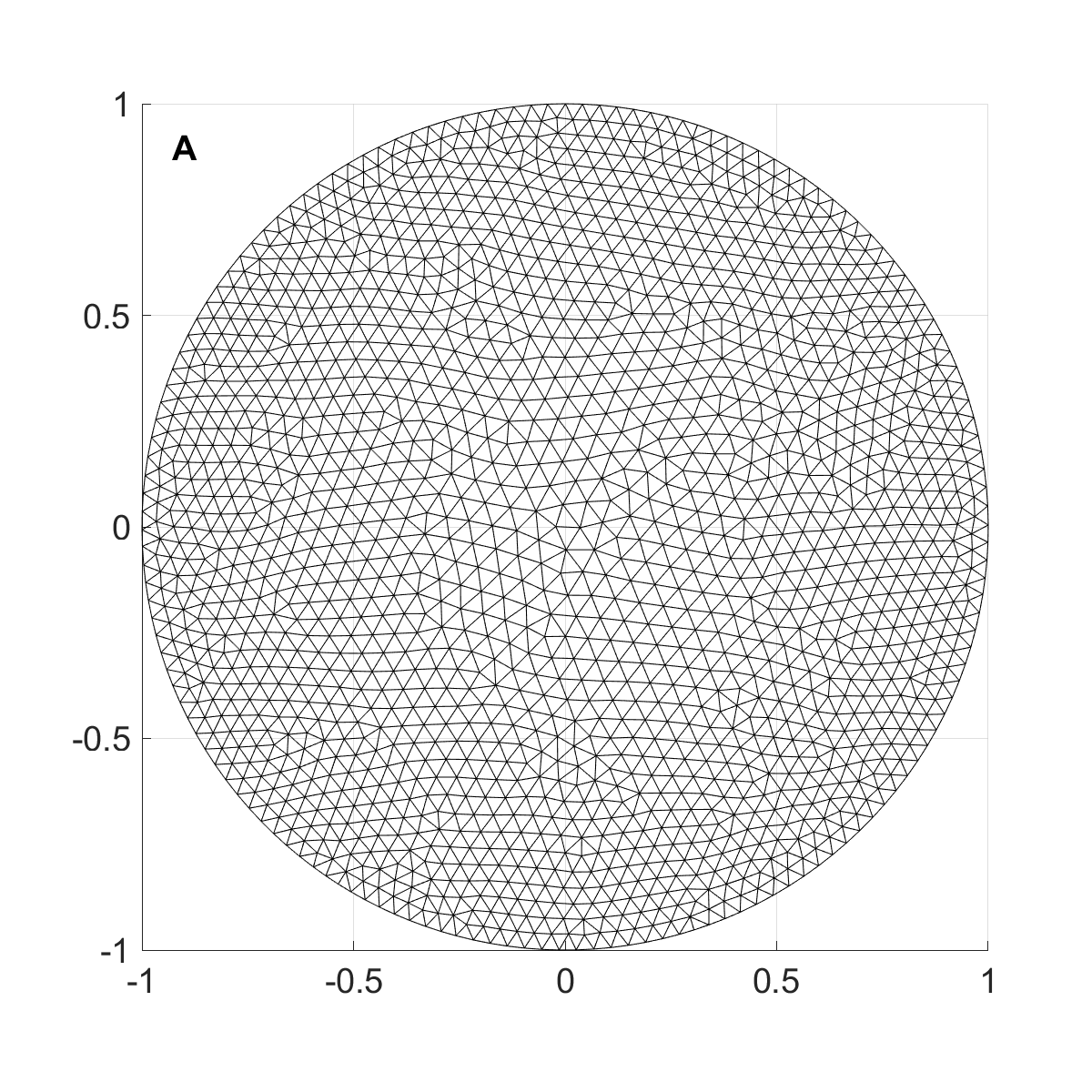}
  \includegraphics[width=0.45\textwidth, trim=10 30 10 45, clip]{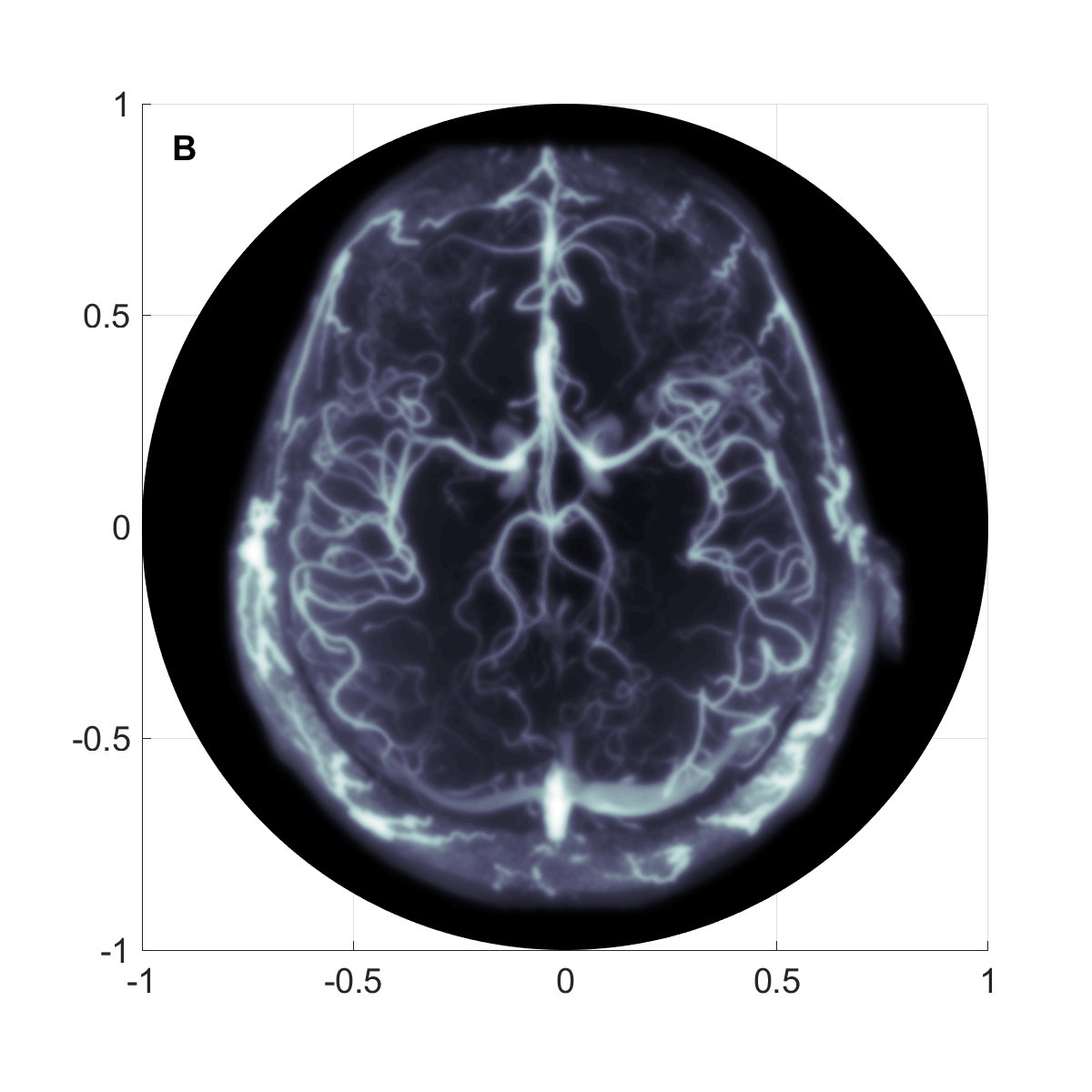}
  \includegraphics[width=0.45\textwidth, trim=10 25 10 45, clip]{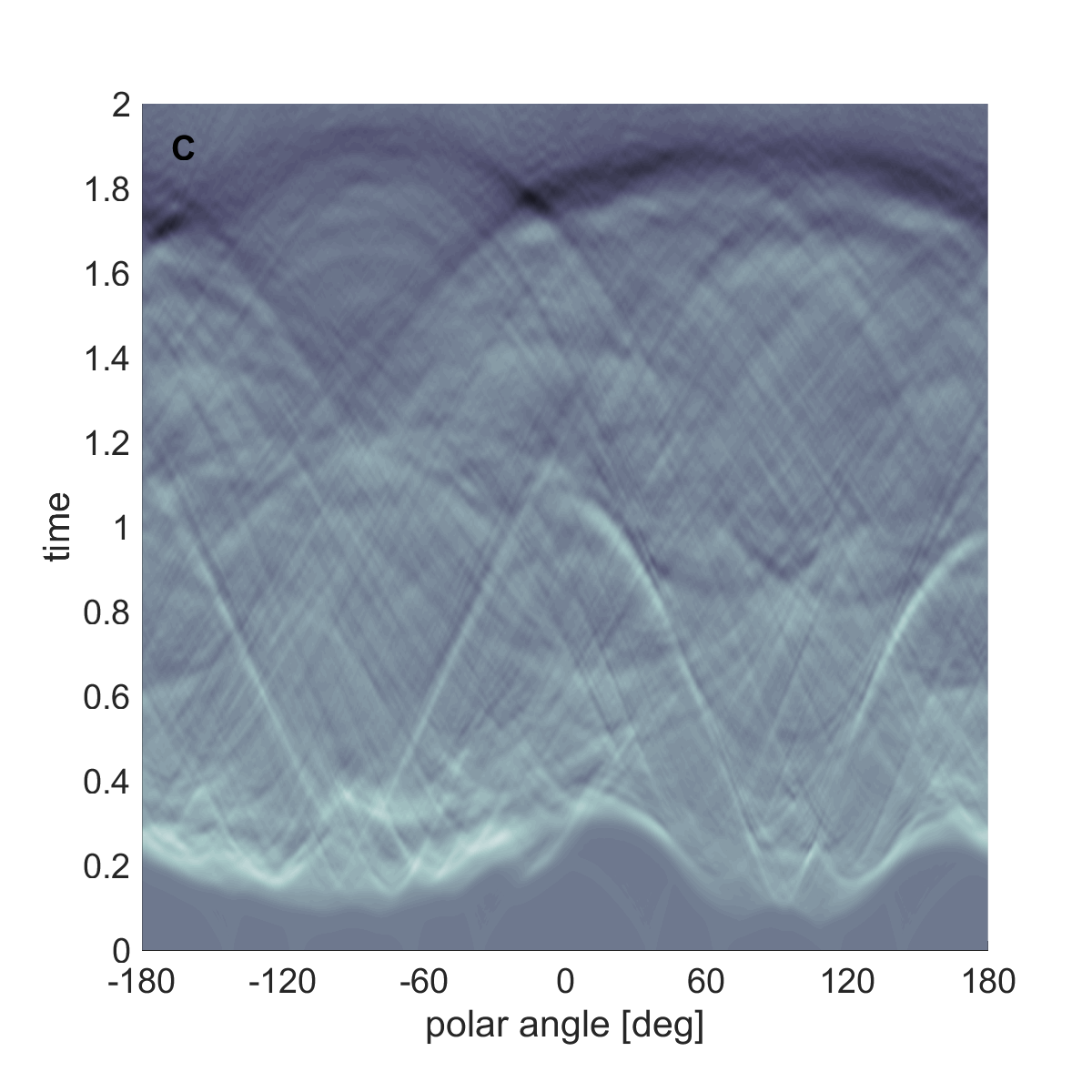}
  \caption{\textit{\textbf{A}: Coarse mesh for the FEM method. \textbf{B}: Exact profile to be reconstructed showing the brain vasculature imaged with MRI technology \cite{Koroshetz2018}. \textbf{C}: Synthetic piezoelectric measurements as modeled by (\ref{eq:forward}).}}
  \label{fig:coarse_mesh}
\end{figure}

Both the forward map $\FF$ and its adjoint $\FF^{*}$ are discretized using a piecewise linear finite element method (FEM) in space and the explicit Newmark method for the time stepping. The discretization parameters are chosen to satisfy the CFL stability condition. The FEM is implemented on triangulations of the domain $\Omega$. For the numerical simulations, we have non-dimensionalized the physical parameters displayed in Table \ref{tab:1} in order to have $c = 1$, $\text{diam}(\Omega) = 2$ and $\rho=1$. The final time $T = 2$. The non-dimensional parameters of the piezoelectric film have  been chosen as follows $\rhosen = 1.5$, $\Csen = 1.0$. The parameters of the backing layer are $\rhosub = 2.0$ and $\Csub = 1.0$. The piezoelectric-elastic coupling coefficient $\kappa = 0.9$. These non-dimensional parameters are consistent with the ranges of their dimensional counterparts listed in Table \ref{tab:1}. 

Figure \ref{fig:coarse_mesh} displays a coarse mesh used for the FEM, the exact pressure profile to be reconstructed and the boundary measurements. These measurements were synthetically generated by applying the discrete version of the forward operator $\FF$ using the aforementioned numerical method for the wave equation. The FEM for the reconstruction procedure has 109,762 degrees of freedom and 6,400 time steps covered the time window for $T=2$. The mesh employed to generate the measurements was more refined, with mesh size approximately half of the mesh size employed in the reconstruction steps, and the data was down-sampled to the reconstruction mesh using linear interpolation.

\begin{figure}[htbp]
  \centering
  \includegraphics[width=0.45\textwidth, trim=0 10 0 30, clip]{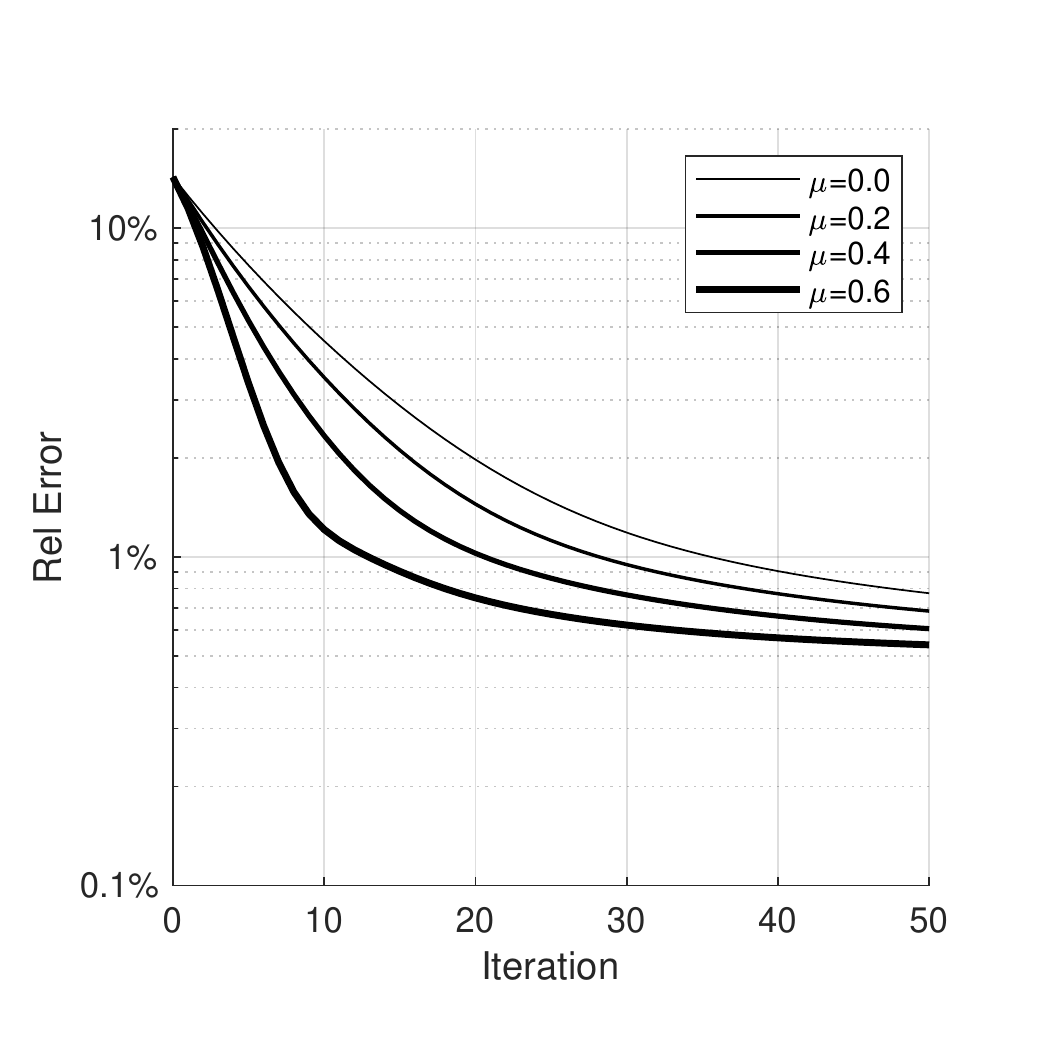}
  \caption{\textit{Relative error versus iteration number, for various values of the momentum factor $\mu$ to accelerate the Landweber iterations. The relaxation factor $\gamma = 5 \times 10^{-2}$ in all cases.}}
  \label{fig:performance_iterations}
\end{figure}

The performance of the Algorithm \ref{alg:landweber} for various values of the momentum factor $\mu$ is displayed in Figure \ref{fig:performance_iterations}. Significant improvements are observed for increasing values of $\mu$. For instance, in order to reach below $1\%$ relative error, the original Landweber method ($\mu = 0.0$) takes 36 iterations, whereas the accelerated method with $\mu=0.6$ takes 12 iterations. For values $\mu > 0.7$, some instability is observed.

\begin{figure}[htb]
  \centering
  \includegraphics[width=0.45\textwidth, trim=20 20 20 20, clip]{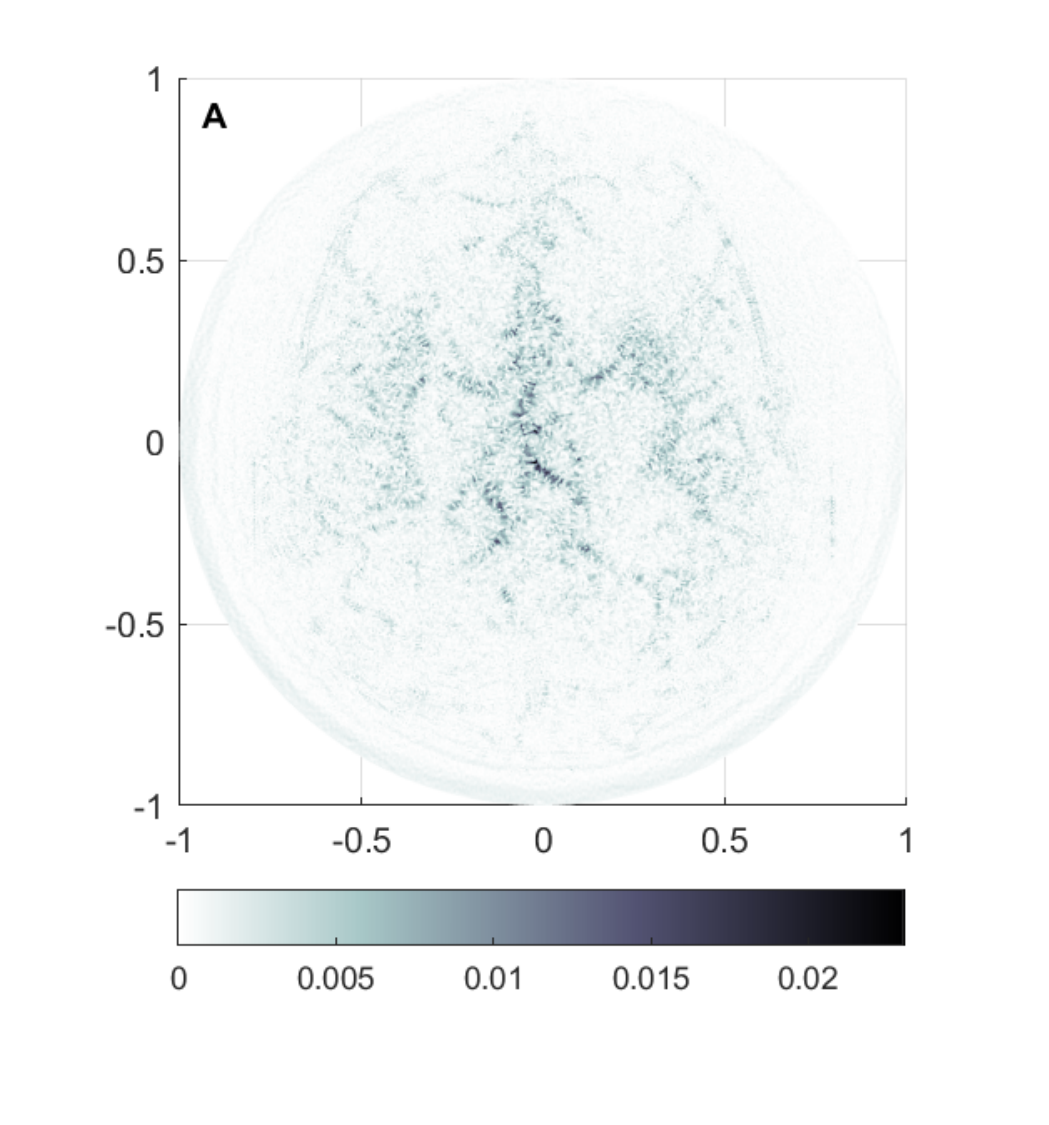}
  \includegraphics[width=0.45\textwidth, trim=20 20 20 20, clip]{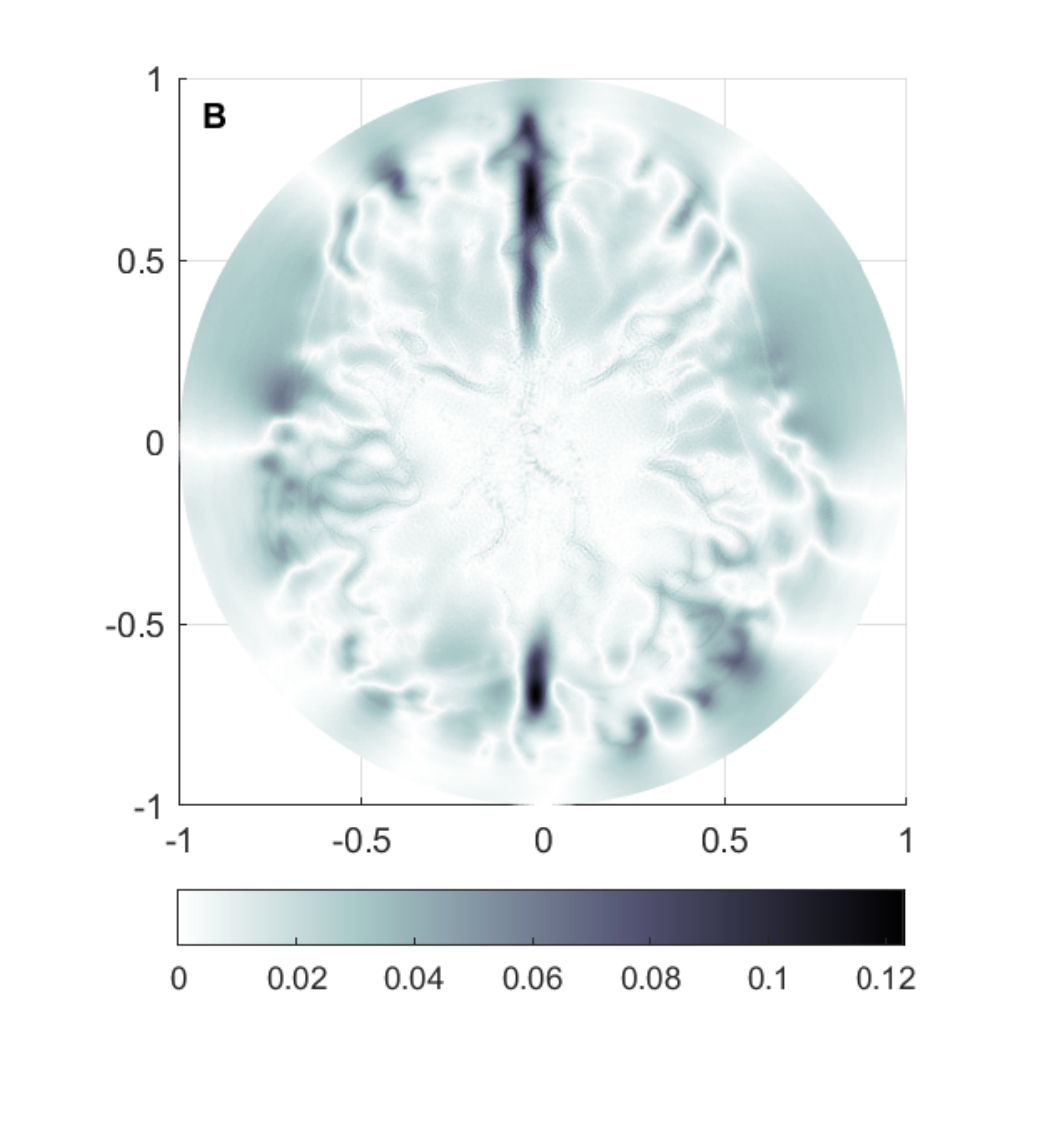}
  \caption{\textit{Error profiles for synthetic measurements obtained using the piezoelectric model (\ref{eq:forward}). \textbf{A}: Reconstruction using an interpretation of the measurements as piezoelectric data. The relative error is $0.54\%$. \textbf{B}: Reconstruction using a naive interpretation of the measurements as Dirichlet data. The relative error is $6.43\%$.}}
  \label{fig:Error_profiles}
\end{figure}

In order to visualize the impact of improperly modeling the piezoelectric measurements, we have implemented two reconstructions of the initial acoustic profile. In both case, the same synthetic measurements are used. For the first reconstruction, we properly interpret the given measurements as generated by the piezoelectric model (\ref{eq:forward}) and carry out the reconstruction using Algorithm \ref{alg:landweber}. After 50 iterations we obtain a relative error of $0.54\%$. For the second reconstruction, we improperly interpret the measurements as the Dirichlet data of the pressure field. This is the naive model commonly employed by others but inconsistent with the piezoelectric transduction. The reconstruction is carried out using Algorithm \ref{alg:landweber} modified by setting $\kappa = 0$ which is equivalent to assuming that the measurements are Dirichlet data. After 50 iterations we obtain a relative error of $6.43\%$. The error profiles for both reconstructions are displayed in Figure \ref{fig:Error_profiles}.

\begin{figure}[htbp]
  \centering
  \includegraphics[width=0.45\textwidth, trim=0 10 0 30, clip]{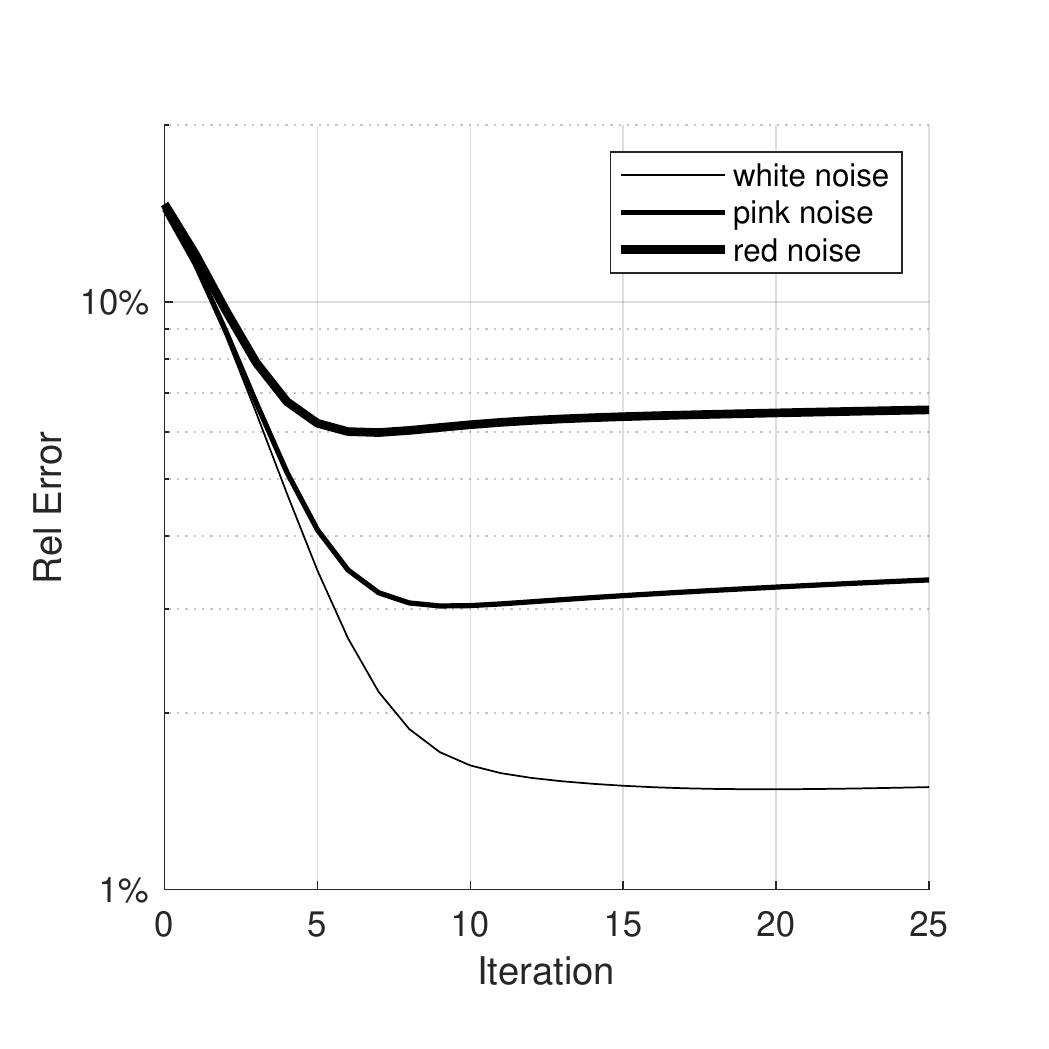}
  \caption{\textit{Relative error versus iteration number for white, pink and red noise. In all cases the noise level is at $10\%$.}}
  \label{fig:noise}
\end{figure} 

The effect of noise added to the piezoelectric measurements is also explored. We have considered three types of noise: white, pink and red. These are characterized by increasing autocorrelation distances or equivalently faster decay of their spectral power. White or Gaussian noise has equal power spectral density (PSD) at different frequencies. Pink noise has a PSD decaying like $\omega^{-1}$ as $\omega \to \infty$. Red or Brownian noise has a PSD decaying like $\omega^{-2}$. Figure \ref{fig:noise} displays the relative error as a function of the iteration number for the three types of noise at a $10\%$ level. We observe that the reconstruction method can handle best the white noise. The error for the pink noise is greater. And the reconstruction error for the red noise is the greatest of the three types of noise.

\begin{figure}[htbp]
  \centering
  \includegraphics[width=0.45\textwidth, trim=0 5 0 10, clip]{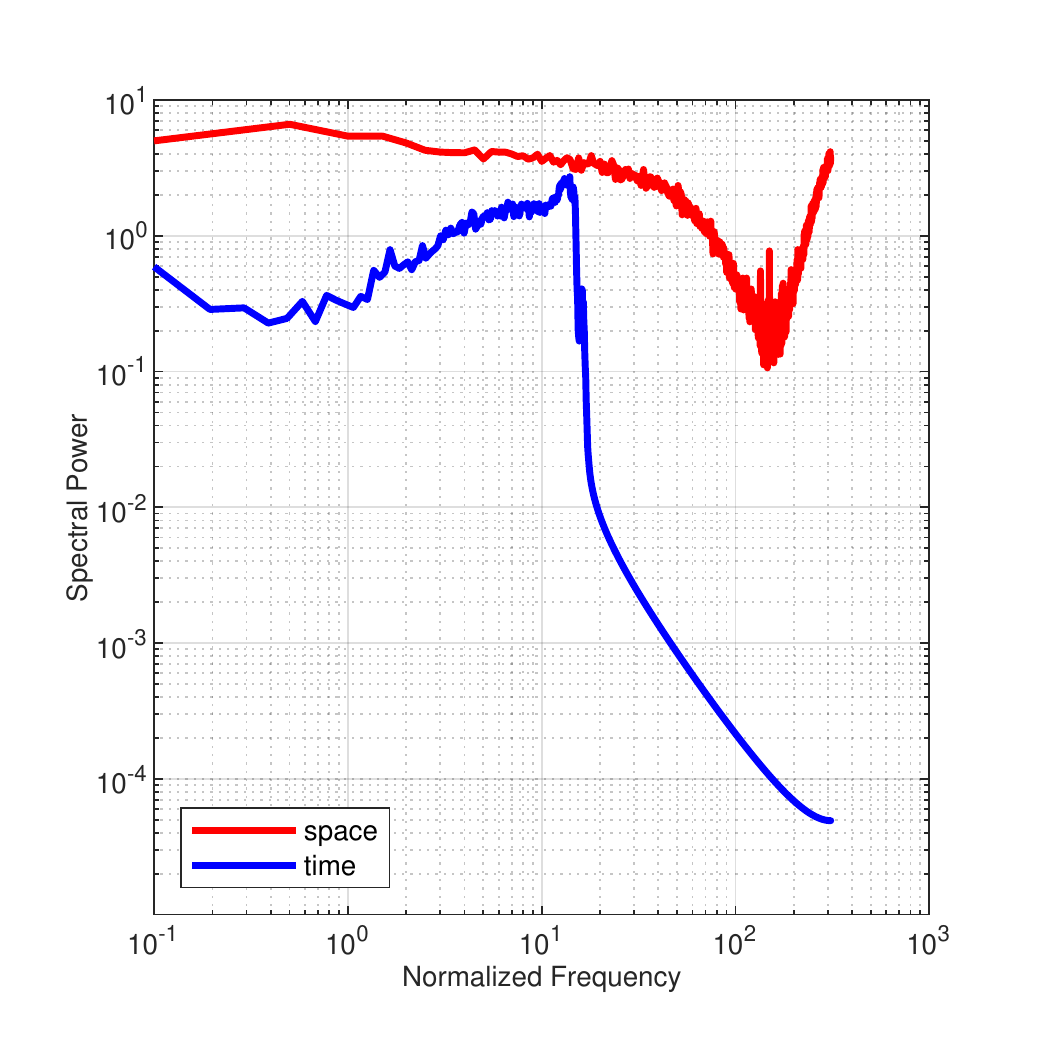}
  \caption{\textit{Frequency response of the normal operator $\left(  \FF \FF^{*} \right)$ to white noise as the input.}}
  \label{fig:noise_power}
\end{figure}

This phenomenon could be explained by studying the spectral characteristics of the discrete normal operators $\left( \FF^{*}\FF \right)$ or $\left( \FF\FF^{*} \right)$. The iterations from the Landweber algorithm correspond to truncated Neumann series \cite{Stefanov2017a},
\begin{equation*}
u_{K} = \gamma \sum_{k=0}^{K} \left( I - \gamma \FF^{*}\FF \right)^{k}(\FF^{*}V),
\end{equation*}
where binomial formula yields
\begin{equation*}
\left( I - \gamma \FF^{*}\FF \right)^{k} = \sum_{j=0}^{k} {k \choose j} \left( - \gamma \FF^{*} \FF \right)^{j}.
\end{equation*}
Therefore, the Landweber iterate $u_{K}$ could be expressed as
\begin{IEEEeqnarray}{ll} \label{eq:NeumannSeries2}
u_{K} &= \sum_{k=0}^{K} \sum_{j=0}^{k} c_{kj} \left( \FF^{*} \FF \right)^{j} (\FF^{*}V)\nonumber \\
 &= \FF^{*} \sum_{k=0}^{K} \sum_{j=0}^{k} c_{kj} \left(  \FF \FF^{*} \right)^{j} V
\end{IEEEeqnarray}
for some coefficients $c_{kj}$. This last expression motivates the study of how the discrete normal operator $\left( \FF \FF^{*} \right)$ responds to the noise contained in the piezoelectric measurements. Figure \ref{fig:noise_power} displays the average power spectral response to time- and space-tracings of white noise as the input to $\left(  \FF \FF^{*} \right)$. We observe that the operator $\left( \FF \FF^{*} \right)$ suppresses the high frequency components. This explains why the reconstruction algorithm can handle white noise better than pink or red noise. The latter noise types have a larger portion of their power residing over low frequencies. Therefore, the reconstruction algorithm does not suppress those types of noise.

\section{Discussion and Conclusions} \label{Section.Conclusions}

We have proposed a mathematical model for the acoustic measurements transduced by piezoelectric sensors. This model, taking the form (\ref{eq:forward}), is highly idealized as described in the Introduction. This idealization allows us to attain a balance between correctness and simplicity in order to analyze the properties of the model (see section \ref{Section.ModelWaves}) and the solvability of the associated photoacoustic problem (see sections \ref{Section.ForwardProblem}-\ref{Section.InverseProblem}). 

The directional response for plane waves derived in section \ref{Section.ModelWaves} and displayed in Figure \ref{fig:directivity} matches the experimental measurements carried out by other researchers \cite{Wilkens2007,Wear2018b,Wear2019a,Brown2000}. Design considerations for the mechanical and electrical properties of the piezoelectric film (condensed in the parameter $\kappa$ defined in (\ref{eq:Kappa})) play a role in the appearance of critical angles where the sensitivity of the sensor vanishes. The incorporation of this type of sensor response into reconstruction algorithms has been highlighted as one of the challenges associated with photoacoustic imaging \cite{Wang2007,Cox2009a,Ellwood2014,Wang2011a,Xia2014} and partially investigated in \cite{Finch2005p,Acosta-Montalto-2015,Zangerl2019,Dreier2019,Acosta2019a}. Our present paper represents a novel contribution to this research effort. 

We note that our mathematical model for the piezoelectric sensor is qualitatively similar to the model for the Fabry--Perot transducer proposed in \cite{Acosta2019a} in spite of the completely different physical principles from which they are derived. Hence, a common mathematical framework can be used to analyze both types of sensors. From the theoretical perspective, as stated in Theorem \ref{Thm:Main}, we conclude that the photoacoustic imaging problem is solvable for piezoelectric measurements. However, some stability is lost compared to measuring directly the Dirichlet data. 

We have also proposed an iterative reconstruction Algorithm \ref{alg:landweber} based on an accelerated Landweber method. We implemented proof-of-concept synthetic simulations to highlight the importance of incorporating the proper modeling of the piezoelectric transducer. We carried out reconstructions with and without the proper piezoelectric model. See the error profiles shown in Figure \ref{fig:Error_profiles}. We observe that certain features, whose wavefronts may have reached the measurement boundary at a non-normal incidence, are missed by the naive reconstruction. We also investigated the effect of noise of different spectral characteristics. The reconstruction method can handle white noise better than pink or red noise. See Figures \ref{fig:noise}. This is due to suppression of high-frequency components in the measurements. See Figure \ref{fig:noise_power}.

One limitation of the proposed method relates to the numerical characteristics of the FEM and Newmark time-stepping. Notice in Figure \ref{fig:performance_iterations} that the error initially decays exponentially as the theory for the Landweber method predicts. However, as the iterations continue, the error eventually stagnates. This stagnation may be attributed to the fact that the measurements were synthetically generated on a mesh different from the reconstruction mesh. Among other issues, the discrete version of $\FF^{*}$ may not be an exact adjoint for the discrete version of $\FF$. Also, the numerical scheme to solve the wave equation suffers from numerical dispersion. Different frequency components of the initial pressure profile travel at different group velocity towards the detection boundary. As a consequence, the discrete version of $\left(\FF^{*}\FF\right)$ loses coercivity (becomes ill-conditioned) and the reconstructed images suffer from aberration. Remedies for this phenomenon have been investigated, including regularization, two-grid methods and numerical schemes with lower dispersion. See \cite[Sect. 6.8-6.10]{GlowinskiLionsHe2008}, \cite{Zuazua2005} and references therein. However, these improvements fall outside of the scope of this paper.

As future research, it would be very important to explicitly model the influence of matching layers commonly employed in the design of piezoelectric sensors \cite{Brown2000,Szabo2014}. Matching layers play an important role in optimizing the transmission of high-frequency acoustic energy into the sensor. In this paper, this transmission is modeled by (\ref{eq:IVBP-02}) which for a flat surface simplifies to
\begin{equation*}
\frac{1}{\rho} \frac{\partial p}{\partial n} + \frac{1}{Z_{\rm b}} \frac{\partial p}{\partial t} = 0.
\end{equation*}
Here $Z_{\rm b} = \rhosub \Csub$ is the acoustic impedance of the thick backing substrate which usually does not match the acoustic impedance $Z = \rho c$ of the fluid or the acoustic impedance $Z_{\rm p} = \rhosen \Csen$ of the piezoelectric film. Matching layers are designed to have an impedance $Z_{\rm ml}$ such that $Z < Z_{\rm ml} < Z_{\rm p} \lesssim  Z_{\rm b}$ so that the wave field experiences a more gradual change of media across the fluid-sensor interface.

The validity of the proposed asymptotic model is limited to small values for the thickness $\epsilon$ of the piezoelectric film with respect to the wavelength of the acoustic waves. Advanced industrial processes are able to manufacture piezoelectric films with thickness in the range 30 -- 100 $\mu$m approximately. As shown in Table \ref{tab:1}, the wave speed in PVDF materials ranges from 1300 -- 2300 m/s. Hence, we expect our model to be valid for frequencies $\lesssim$ 12 MHz. For higher frequencies, the wave field is affected by the presence of the piezoelectric film and the thin matching layer \cite{Brown2000}. In such a scenario, our transmission model would be inappropriate. Hence, there remains a need for a more complete model that can accurately handle multiple frequency scales.

\section*{Acknowledgment}
The author would like to
thank Texas Children's Hospital for its support and for the research oriented environment
provided by the Predictive Analytics Laboratory. The author gratefully acknowledges support by the NSF grant DMS-1712725. The author also thanks the anonymous reviewers whose observations were most constructive.

\ifCLASSOPTIONcaptionsoff
  \newpage
\fi



\bibliographystyle{IEEEtran}
\bibliography{InvProblemBiblio}

%

\begin{IEEEbiography}[{\includegraphics[width=1in,height=1.25in,clip,keepaspectratio]{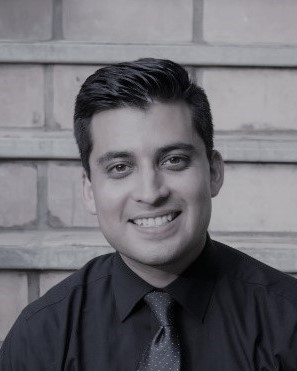}}]{Sebastian Acosta}
received a B.Sc. degree
in mechanical engineering and a M.Sc. degree in
mathematics from Brigham Young University, Provo, Utah, USA in 2009 and 2011 respectively, 
and a Ph.D. degree in computational and applied mathematics from William Marsh Rice University, Houston, Texas, USA in 2014. He is currently a research faculty member in the Department of Pediatrics at Baylor College Medicine and a member of the Predictive Analytics Laboratory at Texas Children's Hospital. His current research interests include biomedical imaging and advanced machine learning for medical monitoring of critically ill patients.
\end{IEEEbiography}







\end{document}